\documentclass[12pt]{article}
\usepackage{amssymb}
\usepackage{graphicx}
\usepackage{amsmath}

\newtheorem{theorem}{Theorem}

\newtheorem{lemma}[theorem]{Lemma}

\newenvironment{proof}[1][Proof]{\noindent\textbf{#1.} }{\QED}

\def\QED{\hfill \ \rule{0.5em}{0.5em}}

\newcommand{\Z}{\mathbb Z}
\DeclareMathOperator{\diam}{\mathrm{diam}}
\DeclareMathOperator{\Cay}{\mathrm{Cay}}
\DeclareMathOperator{\Aut}{\mathrm{Aut}}

\begin{document}

\title{Finite groups of uniform logarithmic diameter}
\author{Mikl\'{o}s Ab\'{e}rt\footnote{Email: 
abert(at)math(dot)uchicago(dot)edu.
Partially supported by NSF Grant DMS-0401006.}
 and L\'{a}szl\'{o} Babai\footnote{
Email:\ laci(at)cs(dot)uchicago(dot)edu.}}
\date{May 20, 2005}
\maketitle

\begin{abstract}
We give an example of an infinite family of finite groups $G_{n}$ such
that each $G_{n}$ can be generated by $2$ elements and the diameter 
of every Cayley graph of $G_{n}$ is $O(\log (\left| G_{n}\right| ))$. 
This answers a question of A.~Lubotzky.
\end{abstract}

Let $G$ be a finite group and $X$ a set of generators. 
Let $\Cay (G,X)$ denote the undirected Cayley graph 
of $G$ with respect to $X$, defined by having vertex set $G$
and $g\in G$ being adjacent to $gx^{\pm 1}$ for $x\in X$.

We define $\diam(G,X)$ to be diameter of $\Cay (G,X)$, that is,
the smallest $k$ such that every element of $G$ can be expressed as a word
of length at most $k$ in $X$ (inversions permitted). 
The diameter of a Cayley graph 
is related to its isoperimetric properties (cf. \cite{aldous, bsz, lubok}).

Let us define the \emph{worst diameter} of $G$, 
\begin{equation*}
\diam_{\max}(G)=\max \left\{ \diam (G,X)\mid X\subseteq
G,\left\langle X\right\rangle =G\right\}
\end{equation*}
to be the maximum diameter over all sets of generators of $G$.

Alex Lubotzky \cite{lubo} asked whether there exists a family of groups $%
G_{n}$ such that the $G_{n}$ have a bounded number of generators and
the worst diameter of $G_{n}$ is $O(\log \left|
G_{n}\right| )$. 

Lubotzky's question has been addressed by Oren Dinai \cite
{dinai} who proved that the family $G_{n}=SL_{2}(\Z/p^{n}\Z)$ (where $p$ is a
fixed prime) has worst diameter which is poly-logarithmic in the size of $%
G_{n}$. Her result relies on the work of Gamburd and Shahshahani \cite
{gamburd} who proved the corresponding result in the case when the
generating set projects onto $G_{2}$.

We demonstrate that there exists a family of solvable groups with
logarithmic worst diameter.   Let $p$ be an odd prime and let $W(p)$
denote the wreath product $W_{p}=C_{2}\wr C_{p}$. %% of order $p2^p$.
%%$C_{2}^p$ by $C_p$).

\begin{theorem}  \label{thm:main}
\qquad  $\diam_{\max}(W_{p}) \le 20(p-1)$.
\end{theorem}

Since the groups $W_{p}$ can be generated by $2$ elements and the size of $%
W_{p}$ is $\left| W_{p}\right| =p2^{p}$, this answers Lubotzky's question
affirmatively.

\bigskip

Let us fix some notation.  Let $F=\mathbb{F}_{2}$ denote the field of $2$
elements.  The wreath product $W=W_{p}$ is the
semidirect product of the $F$-vector-space $U\vartriangleleft W$ of dimension 
$p$ and the cyclic group $C=C_{p}<W$ of order $p$. The structure of the 
$C$-module $U$ is governed by the irreducible factors of the polynomial 
$x^{p}-1$ over $F$.  It turns out that $U$ splits as 
\begin{equation*}
U=T\times M_{1}\times \cdots \times M_{k}
\end{equation*}
where $T$ is a trivial (one-dimensional) module and the $M_{i}$ are
nontrivial, pairwise inequivalent simple modules. This follows from the
fact that $x^{p}-1$ has no multiple roots over $F$. Since the Galois
group of the corresponding extension of $F$ is generated by the Frobenius
automorphism $x\longmapsto x^{2}$, the dimension of each simple module $%
M_{i}$ is equal to the multiplicative order of $2$ modulo $p$.

The center of $W$ is $T$. It will be more convenient to factor
out $T$ and compute the worst diameter in the quotient group $W/T$ first.

Let $G=W/T$. Then $G$ is the semidirect product of the $F$-vectorspace $%
V\vartriangleleft W$ of dimension $p-1$ and $C$. Let us
fix a generator $c$ of $C$. We shall write the elements of $G$ in the form 
$vc^{i} $ where $v\in V$ and $i\in \{0,\dots ,p-1\}$. The $C$-module $V$
splits as 
\begin{equation*}
V=M_{1}\times \cdots \times M_{k}
\end{equation*}
where the $M_{i}$ are as above.

\bigskip

\begin{lemma}
\label{tran}The automorphism group $\Aut (G)$ acts transitively on $%
G\backslash V$.
\end{lemma}

\begin{proof}
Using the identity 
\begin{equation*}
v^{-1}cv=v^{-1}(0c)v=(v^{c^{-1}}-v)c=v^{c^{-1}-1}c,
\end{equation*}
we see that the element $c$ is conjugate to an arbitrary element $wc$ where $%
w\in V^{(c^{-1}-1)}=V$ (here we use that the $M_{i}$ are simple nontrivial
modules).

Now the wreath product $W_{p}$ can be understood as the semidirect product
of the group algebra $FC_{p}$ by $C_{p}$. This shows that every
automorphism of $C_{p}$ extends to an automorphism of $W$.  Since the
center is characteristic, the same holds for $G$. That is, the $p-1$
conjugacy classes outside $V$ collapse into one automorphism class. 
\end{proof}

\bigskip
The following well-known observation appears, e.\,g., as
\cite[Lemma 5.1]{bs}.

\begin{lemma}
\label{reszcsop}Let\thinspace\ $G$ be a finite group and let $%
N\vartriangleleft G$. Then 
\begin{equation*}
\diam_{\max}(G)\leq 2\diam_{\max}(G/N)\diam_{\max}(N)+\diam_{\max}%
(N)+\diam_{\max}(G/N)
\end{equation*}
\end{lemma}

\begin{proof}
For completeness, we include the proof.  
Let $X$ be a set of generators of $G$. Then there exists a set $T$ of 
coset representatives of $N$ in $G$ such that every element of $T$ 
can be expressed as a word of length at most 
$\diam_{\max}(G/N)$ in $X$. Now the set $S$ 
of Schreier generators, defined as
\begin{equation*}
S=\left\{ txu^{-1}\mid t\in T,x\in X,u\in txN\cap T\right\},
\end{equation*}
generates $N$, so every element of $N$ can be expressed as a word of length 
\begin{equation*}
  \le(2\diam_{\max}(G/N)+1)\diam_{\max}(N)
\end{equation*}
in $X$. Finally, $G=TN$ gives us the required estimate.
\end{proof}

\bigskip

We note that the lemma also holds if the subgroup $N\leq G$ is not normal; 
in this case $\diam_{\max}(G/N)$ should denote the worst diameter of a 
Schreier graph of $G$ with stabilizer $N$.

\bigskip

\begin{theorem} \label{igazi}
\qquad  $\diam_{\max}(G_{p})\leq \frac{13}{2}(p-1)$.
\end{theorem}

The essence of the proof will be contained in the 
following case which refers to 
a specific type of generating set.

\begin{lemma}   \label{lemma:special}
Let $X=\{c, w_2,\dots,w_n\}$ where $w_2,\dots,w_n\in V$ and assume 
$X$ generates $G_p$.  Then every $v\in V$ can be obtained as
a word of length $\le 3(p-1)$ in $X$ with the $w_i$ occurring at
most a total of $p-1$ times.
\end{lemma}

% ****************

\begin{proof}
We may assume that $c,w_{2},\dots ,w_{n}$ is an irredundant generating
set; otherwise, drop some of the $w_i$. For $2\leq i\leq n$ let 
\begin{equation*}
A_{i}=\left\{ j\mid w_{i,j}\neq 0\right\}
\end{equation*}
and 
\begin{equation*}
B_{i}=A_{i}\backslash \bigcup_{j=2}^{i-1}A_{j}.
\end{equation*}
Then the $B_{i}$ are disjoint, non-empty (because of the irredundancy) and 
\begin{equation*}
\bigcup_{j=2}^{n}B_{j}=\left\{ 2,\dots, n\right\}
\end{equation*}
Also let 
\begin{equation*}
V_{j}=\prod_{i\in B_{j}}M_{i}\subseteq M_{1}\times \cdots \times M_{k}=V.
\end{equation*}
Then 
\begin{equation*}
V=V_{2}\times \cdots \times V_{n}.
\end{equation*}
Let $d_{j}=\dim _{F}V_{j}$.

Let $F[x]_{d}$ denote the space of polynomials of degree at
most $d$ over $F$. Since $V_{j}$ is a direct product of simple pairwise
non-equivalent $F[x]$-modules, if $v\in V_{j}$ is a generator of $V_{j}$
then 
\begin{equation*}
V_{j}=vF[x]=vF[x]_{d_{j}-1}
\end{equation*}
(otherwise a polynomial of degree at most $d_{j}-1$ would annihilate
$V_{j}$, a contradiction).

This implies that 
\begin{equation*}
V=w_{2}F[x]_{d_{2}-1}+\cdots +w_{n}F[x]_{d_{n}-1}.
\end{equation*}

Changing to group notation this means that for all $v\in V$ there exist
polynomials $f_{2},\dots ,f_{n}$ such that $\deg f_{j}\leq d_{j}-1$ and 
\begin{equation*}
v=w_{2}^{f_{2}(c)}+\cdots +w_{n}^{f_{n}(c)}
\end{equation*}

Now we will use the Horner scheme to obtain $v$ as a short word in
the generating set $c,w_{2},\dots ,w_{n}$.

We claim that if $w\in V$ and $f(x)\in F[x]$ of degree $d-1$ then $w^{f(c)}$
can be obtained as a word in $c$ and $w$ of length at most $3d$ and we use $%
w $ at most $d$ times. This goes by induction on $d$. For $d=1$ it is
obvious. If $d>1$ then $f(c)=cg(c)+\epsilon $ where $\epsilon =0$ or $1$ and 
$\deg g=d-2$. Now 
\begin{equation*}
w^{f(c)}=c^{-1}w^{g(c)}c+\epsilon w
\end{equation*}
which by induction has length at most $2+3(d-1)+1=3d$ (and we used $w$ at
most $d$ times). This proves the claim.

In particular $w_{j}^{f_{j}(c)}$ can be obtained as a word in $c$ and $w_{j}$
of length at most $3d_{j}$. Adding up, $v$ can be obtained as a word in $%
c,w_{2},\dots ,w_{n}$ of length at most $3(p-1)$ where we used the $w_{j}$
at most $p-1$ times. 
\end{proof}

\bigskip
Now we turn to the proof of Theorem~\ref{igazi}.

\bigskip
\noindent \textbf{Proof of Theorem~\ref{igazi}.}
Let $v_{1}c_{1},\dots ,v_{n}c_{n}$ be a set of generators of $G$. 
For $\alpha \in \Aut(G)$ the Cayley graphs  
\begin{equation*}
\Cay (G,\{(v_{1}c_{1})^{\alpha },\dots ,(v_{n}c_{n})^{\alpha }\})%
\text{ and }\Cay (G,\{v_{1}c_{1},\dots ,v_{n}c_{n}\})
\end{equation*}
are isomorphic. Since at least one of the $c_{i}$ have to be nontrivial,
using Lemma~\ref{tran} we can assume that $v_{1}=0$ and $c_{1}=c\neq 1$. 

Now $c,v_{2}c_{2},\dots ,v_{n}c_{n}$ generate $G$ if and only if $%
c,v_{2},\dots ,v_{n}$ do.   

Let 
\begin{equation*}
v_{i}=(v_{i,1},\dots ,v_{i,k})
\end{equation*}
be the decomposition of $v_i$ 
according to $V=M_{1}\times \cdots \times M_{k}$. It is
easy to see that $c,v_{2},\dots ,v_{n}$ generate $G$ if and only if $%
v_{2},\dots ,v_{n}$ generate $V$ as a $C_{p}$-module and this happens if
and only if for all $j$ ($1\leq j\leq k$) there exists $i$ such
that $v_{i,j}\neq 0 $.

For $2\leq i\leq n$ let us define 
\begin{equation*}
w_{i}=[c,v_{i}c_{i}]=
c^{-1}c_{i}^{-1}v_{i}cv_{i}c_{i}=v_{i}^{c_{i}c-c_{i}}=(v_{i}^{c_{i}})^{c-1}.
\end{equation*}
Since the $M_{j}$ are nontrivial simple, $w_{i,j}=0$ if and only if $%
v_{i,j}=0$. Using the observation made in the previous paragraph,
this shows that $c,w_{2},\dots ,w_{n}$ also generate $G$.

Let us now apply Lemma~\ref{lemma:special} to this
latter set of generators.  Noting that 
$w_{j}=[c,v_{j}c_{j}]$ can be obtained as a word of
length $4$ in $c$ and $v_{j}c_{j}$,
we infer from Lemma~\ref{lemma:special} that any $v\in V$ can be obtained as a 
word in the original generating set $%
v_{1}c_{1},\dots ,v_{n}c_{n}$ of length at most $6(p-1)$. So the diameter
of $G$ with respect to $v_{1}c_{1},\dots ,v_{n}c_{n}$ is at most 
\begin{equation*}
6(p-1)+\frac{p-1}{2}=\frac{13}{2}(p-1)\text{.} 
\end{equation*}
\QED

\bigskip

\noindent \textbf{Proof of Theorem~\ref{thm:main}.} The center $T=Z(W)$ has
order\thinspace $2$ so $\diam_{\max}(T)=1$. Using Theorem~\ref{igazi} and
Lemma \ref{reszcsop} we get 
\begin{equation*}
\diam_{\max}(W)\leq 3\diam_{\max}(G)+1\leq \frac{39}{2}(p-1)+1
  \le  20(p-1)
\end{equation*}
what we wanted to show.  \hfill  \QED

\bigskip

\noindent \textbf{Remark}. A similar proof shows that the wreath product $%
C_{q}\wr C_{p}$, where $q$ is a fixed prime and $p$ runs through all
primes not equal to $p$, has worst diameter 
\begin{equation*}
\diam_{\max}(C_{q}\wr C_{p})\leq K_q p
\end{equation*}
where $K_q$ depends only on $q$. 

\bigskip 

For a finite group $G$ let $r(G)$ denote the largest size of irredundant
generating sets.  (A generating set $X$ is irredundant if no proper
subset of $X$ generates $G$). This measure has been investigated by Saxl and
Whiston (see \cite{saxwhi} and references therein) for various classes of
groups.

It is natural to ask the value of $r(G_{p})$.  As we have seen, generation
in $G_{p}$ is governed by the structure of the underlying module and so 
$r(G_{p})=1+k$ where $k=(p-1)/o_{p}(2)$, where $o_{p}(2)$ denotes
the multiplicative order of $2$ modulo $p$.

This permits a wide range of values for $r(G_{p})$ in terms of $p$:
\begin{equation}
           2\le  r(G_{p}) \le 1+\frac{p-1}{\log_{2}(p-1)}.
\end{equation}
Let us look at the extremes.

The minimum possible value $r(G_{p})=2$ occurs if and only if $2$
is a primitive root modulo $p$. Heath-Brown's solution \cite{heath}
to Artin's conjecture tells us that, with maybe two exceptions, every 
%%nonsquare natural number
prime $q$ is a primitive root modulo $p$ for infinitely 
many primes $p$. So with any luck (that is, if $2$ is not one of 
these exceptions), we have
infinitely many primes such that $r(G_{p})=2$. In case we are not lucky, one
of $C_{3}\wr C_{p}$ and $C_{5}\wr C_{p}$ will do, according to
the remark above. 

Conversely, the smallest possible value for $o_{p}(2)$ is 
$\log _{2}(p-1)$ and this occurs if and only if $p$ is
a Mersenne prime.   In this case, $r(G_{p})$ takes its largest 
possible value, $(p-1)/\log _{2}(p-1)$.

%% The best known estimate in this
%% direction [**] gives us that there exist infinitely many primes such 
%% that $r(G_{p})\geq ***$.


\bigskip

\begin{thebibliography}{9}
\bibitem{aldous} {\sc D. Aldous:}
    On the Markov chain simulation method for uniform combinatorial
    distributions and simulated annealing.
    {\sl Probability in Engineering and Informational Sciences}
    {\bf 1} (1987), 33-46.

\bibitem{bs}  {\sc L. Babai, \'A. Seress:}
   On the diameter of permutation groups.
       {\sl Europ. J. Comb.} {\bf 13} (1992), 231-243.

\bibitem{bsz}  {\sc L. Babai, M. Szegedy:}
        Local expansion of symmetrical graphs.
      {\sl Combinatorics, Probability, and Computing} {\bf 1} (1992), 1--11.

\bibitem{dinai}  {\sc O. Dinai:} Poly-log diameter bounds for some 
    families of finite groups. Manuscript, 2004.

\bibitem{gamburd} {\sc A. Gamburd and M. Shahshahani:} 
  Uniform diameter bounds for some families of Cayley graphs.
  {\sl Int. Math. Res. Not.} {\bf 71} (2004), 3813--3824.

\bibitem{heath} {\sc D.\,R. Heath-Brown:} 
   Artin's conjecture for primitive roots.
   {\sl Quart. J. Math. Oxford} {\bf 37} (1986), 27-38.

\bibitem{lubok}  {\sc A. Lubotzky:}
   {\sl Discrete groups, expanding graphs and invariant measures.}
  % With an appendix by Jonathan D. Rogawski, 
    Progress in Mathematics 125, Birkh\"{a}user Verlag, Basel, 1994.

\bibitem{lubo}  {\sc A. Lubotzky:} Collected Problems at the Conference on
Automorphic Forms, Group Theory and Graph Expansion. IPAM 2004

%% whinston - saxl?
\bibitem{saxwhi}  {\sc J. Saxl and J. Whiston:} On the maximal size of 
   independent generating sets of $\mathrm{PSL}_{2}(q)$.
    {\sl  J. Algebra} {\bf 258/2} (2002), 651--657.
\end{thebibliography}
\end{document}